\documentclass[11pt,a4paper]{article}
\usepackage{amsmath,amssymb,url}
\usepackage{epsfig,graphicx}
\addtocounter{MaxMatrixCols}{1} 

\newcommand{\Supp}{\mbox{Supp}}

\newtheorem{theorem}{Theorem}
\newtheorem{lemma}[theorem]{Lemma}
\newtheorem{proposition}[theorem]{Proposition}

\newtheorem{definition}[theorem]{Definition}

\newtheorem{example}{Example}
\begin{document}
\begin{center}
{\LARGE{Roots and coefficients of multivariate polynomials\\

 over finite
  fields}}
\end{center}

 \begin{center}{Olav Geil}\\
Department of Mathematical Sciences\\
Aalborg University\\
{\url{olav@math.aau.dk}}
\end{center}
  
\noindent {\it{{\bf{Abstract:}} 
Kopparty and Wang studied in \cite{kopparty2014roots} the relation
between the roots of a univariate polynomial over ${\mathbb{F}}_q$ and
the zero-nonzero pattern of its coefficients. We generalize their
results to polynomials in more variables.}}\\

\section{Introduction}
In~\cite{kopparty2014roots} Kopparty and Wang considered the zero-nonzero
pattern of a univariate polynomial $P(X)$ over ${\mathbb{F}}_q$ and its
relation to the number of roots in ${\mathbb{F}}_q^\ast$. Their main
theorem~\cite[Th.\ 1]{kopparty2014roots} states that a polynomial with
many zeros cannot have long sequences of consecutive coefficients all
being equal to zero. Then in~\cite[Th.\ 2]{kopparty2014roots} they 
gave necessary and sufficient conditions for a product of pairwise
different linear factors to have sequences of zero
coefficients of maximal possible length for any polynomial with prescribed number of roots. In
this note we generalize the abovementioned results to polynomials in
more variables.\\
In Section~\ref{sec2} we start by recalling 
the results by Kopparty
and Wang. In Section~\ref{sec3} we then present and prove the generalizations.

\section{Univariate polynomials}\label{sec2}
The main theorem in~\cite{kopparty2014roots} is their Theorem 1 which
we present in a slightly stronger version.

\begin{theorem}\label{the1}
Let $P(X) \in {\mathbb{F}}_q[X]$ be a nonzero polynomial of degree at
most $q-2$, say $P(X)=\sum_{i=0}^{q-2}b_iX^i$. Let $m$ be the number
of $x\in {\mathbb{F}}_q^\ast$ with $P(x) \neq 0$. Then there does not
exist any $k \in \{0, \ldots , q-2\}$ where all the $m$ coefficients
$b_k$, $b_{k+1 {\mbox{ mod }} (q-1)}, \ldots, b_{k+m-1 {\mbox{ mod }}
  (q-1)}$ are zero.
\end{theorem}

The modification made in Theorem~\ref{the1} is that we consider $k\in
\{0, \ldots , q-2\}$ rather than just $k \in \{0, \ldots ,
q-1-m\}$. The proof in~\cite{kopparty2014roots} is easily modified to
cover this more general situation. Alternatively, one can deduce it by
writing $P(X)=X^sQ(X)$ with $s$ maximal and then
applying~\cite[Th.\ 1]{kopparty2014roots} to $Q(X)$.\\

Obviously, if we consider a product of $q-1-m$ pairwise different
linear factors $X-x$ with $x\neq 0$, this polynomial has exactly $m$
non-roots in ${\mathbb{F}}_q^\ast$ and we have $b_{q-m}=\cdots =b_{q-2}=0$ which is a
sequence of $m-1$ consecutive zero coefficients modulo $q-1$. The
below theorem, corresponding to~\cite[Th.\ 2]{kopparty2014roots}, 
gives sufficient and necessary conditions for 
a sub-sequence of $m-1$ consecutive zeros among $b_0, \ldots ,
b_{q-m-2}$ to exist.
\begin{theorem}\label{the2}
Let $S$ be a subset of ${\mathbb{F}}_q^\ast$ of size $q-1-m$, where $m
\geq 2$ and consider
\begin{equation}
P(X)=\prod_{a \in S}(X-a)=\sum_{i=0}^{q-1-m}b_iX^i. \label{eqtrekant}
\end{equation}
There exists a $k\in \{1, \ldots  , q-2m\}$ such that $b_k= \cdots
=b_{k+m-2}=0$ if and only if ${\mathbb{F}}_q^\ast \backslash S$ is
contained in $\gamma H$ for some $\gamma \in {\mathbb{F}}_q^\ast$ and
for some proper multiplicative subgroup $H$ of ${\mathbb{F}}_q^\ast$. 
\end{theorem}

Inspecting the proof in~\cite{kopparty2014roots} one sees that for
polynomials of the form~(\ref{eqtrekant}) the
existence of one sub-sequence of $m-1$ consecutive zero coefficients in
$b_0, \ldots , b_{|S|-1}$ 
is equivalent to the existence of $(q-1)/|H|$ such disjoint sequences.  

\begin{proposition}\label{pro1}
Let $P(X)$ be a polynomial as in~(\ref{eqtrekant}) satisfying the
condition of Theorem~\ref{the2}. That is, there exists a $k \in \{1,
\ldots , q-2m\}$ such that $b_k=\cdots
=b_{k+m-2}=0$ where $m=| {\mathbb{F}}_q^\ast \backslash S|$. Write
$d=|H|$ where $H$ is the subgroup corresponding to $P$. The
coefficients $b_{jd}, b_{jd+(d-m)}$, $j=0, \ldots , \frac{q-1}{d}-1$ are nonzero and the only other
possible nonzero coefficients of $P(X)$ are $b_{jd+1},
b_{jd+1},\ldots , b_{jd+(d-m)-1}$, $j=0, \ldots ,
\frac{q-1}{d}-1$. 
\end{proposition}

\noindent {\bf{Proof:}} According to~\cite[Proof of Th.\ 2]{kopparty2014roots}, if $P(X)$
satisfies the conditions in Theorem~\ref{the2} then it can be
written 
$$\bigg( \sum_{j=1}^{(q-1)/d}b_jX^{(q-1)-jd}\bigg) \cdot U(X)$$
where $U$ is a product of $d-m$ pairwise different expressions $X-x$
with $x \in {\mathbb{F}}_q^\ast$.

\begin{example}\label{ex1}
Let $\alpha$ be a primitive element of ${\mathbb{F}}_{16}$. We first
consider
$$T=\{\beta \mid \beta^3=\alpha^3\}=\{\alpha,\alpha^6,\alpha^{11}\}.$$
The support of $P(X)$ becomes $\{1, X^3, X^6,X^9,X^{12}\}$. If we
choose $T$ to be a subset of $\{ \alpha, \alpha^6, \alpha^{11}\}$ of
size $2$ then the support of $P(X)$ becomes $\{1,
X,X^3,X^4,X^6,X^7,X^9,X^{10},X^{12},X^{13}\}$. Consider next 
$$T=\{ \beta \mid \beta^5=\alpha^{10} \}=\{\alpha^2,
\alpha^5,\alpha^8,\alpha^{11}, \alpha^{14}\}.$$
The support of $P(X)$ becomes $\{1, X^5,X^{10}\}$. Finally, if we
choose $T$ to be a subset of $\{\alpha^2,
\alpha^5,\alpha^8,\alpha^{11}, \alpha^{14}\}$ of size 3 then we can
conclude:
\begin{eqnarray}
\{1, X^2,X^5, X^7,X^{10},X^{12}\} \subseteq \Supp P \subseteq \{1, X,X^2,X^5, X^6,X^7,X^{10},X^{11}, X^{12}\}.\nonumber
\end{eqnarray}
\end{example}

\section{Multivariate polynomials}\label{sec3}
The crucial observation used in the proof of Theorem~\ref{the1} is
that a univariate polynomial $F(X)$ can at most have $\deg F$
roots. For multivariate polynomials over general fields there does not
exist a similar result as typically such polynomials have infinitely
many roots when the field under consideration is infinite. For
multivariate polynomials over finite fields, however, we do have a
counterpart to the bound used in the proof of Theorem~\ref{the1}. We
describe this bound in terms of roots from $({\mathbb{F}}_q^\ast)^n$ in
Proposition~\ref{pro1} below. To motivate the bound we need a few
results from Gr\"{o}bner basis theory.\\

Let ${\mathbb{F}}$ be a field and $I \subseteq {\mathbb{F}}[X_1, \ldots , X_n]$ an
ideal. Throughout this section assume that an arbitrary fixed monomial
ordering $\prec$ has been chosen. Following~\cite{onorin} we define the
footprint of $I$ by
\begin{eqnarray}
\Delta_\prec (I) &=&\{ X_1^{i_1} \cdots X_n^{i_n} \mid X_1^{i_1}
\cdots X_n^{i_n} {\mbox{ is not }} \nonumber \\
&&\, \, \, \, \, \, \, \, {\mbox{ a leading monomial of any polynomial
    in }}I\}.\nonumber
\end{eqnarray}
From~\cite[Prop.\ 4, page 229]{clo} we know that $\{M+I \mid M \in
\Delta_{\prec}(I)\}$ constitutes a basis for ${\mathbb{F}}[X_1, \ldots , X_n]/I$
as a vector space over ${\mathbb{F}}$. Assume $I$ is finite dimensional (which
simply means that $\Delta_\prec(I)$ is a finite set). Consider $\ell$
pairwise different points $P_1, \ldots , P_\ell$ in the zero-set of
$I$ (over ${\mathbb{F}}$). The map ${\mbox{ev}}: {\mathbb{F}}[X_1, \ldots , X_n]/I
\rightarrow {\mathbb{F}}^\ell$ given by ${\mbox{ev}}(F+I)=(F(P_1), \ldots ,
F(P_\ell))$ is a surjective vector space homomorphism (surjectivity
follows by Lagrange interpolation). Therefore 
\begin{equation}
\ell \leq
|\Delta_\prec(I)|\label{eqabove}
\end{equation} 
(this result is often called the footprint bound
\cite{onorin}). In particular we derive:
\begin{proposition}\label{pro2}
Consider $P(\vec{X}) \in {\mathbb{F}}_q[X_1, \ldots , X_n]$ with
leading monomial equal to $X_1^{i_1}\cdots X_n^{i_n}$ such that $i_s
< q-1$ for $s=1, \ldots ,n$. Let $m$ be the number of elements in
$({\mathbb{F}}_q^\ast)^n$ which are not roots of $P$. Then $m \geq
\prod_{s=1}^n (q-1-i_s)$.
\end{proposition}

\noindent {\bf{Proof:}} The proor follows by applying (\ref{eqabove}) to the ideal $I=\langle
P(\vec{X}), X_1^{q-1}-1, \ldots , X_n^{q-1}-1 \rangle$. The footprint
of this ideal 
is a subset of 
\begin{eqnarray}
\{X_1^{j_1} \cdots X_n^{j_n} \mid 0 \leq j_s < q-1,
s=1, \ldots , n,  {\mbox{ not all }} j_s {\mbox{ satisfy }} i_s \leq j_s\}.\nonumber
\end{eqnarray}
Therefore,
the number of non-roots is at least $| \{ (j_1, \ldots , j_n) \mid i_s
\leq j_s < q-1, s=1, \ldots , n\}|$.

Observe that for $n=1$ the statement in Proposition~\ref{pro2} is
but the well-known fact that a multivariate polynomial $P$ has
at least $q-1-\deg P$ non-roots in ${\mathbb{F}}_q^\ast$.\\

Before giving the generalization of Theorem~\ref{the1} we introduce
the set $U(q,m,n)$. This set shall play the role as did the set of consecutive monomials
$\{X^{q-1-m}, \ldots , X^{q-2}\}$ in connection with Theorem~\ref{the1}.

\begin{definition}
Given positive integers $m$ and $n$ let
$${\mathcal{M}}(q,n)=\{ X_1^{i_1}\cdots X_n^{i_n} \mid 0 \leq i_1,
\ldots , i_n < q-1\},$$
 $$U(q,m,n)=\{X_1^{i_1} \cdots X_n^{i_n} \in {\mathcal{M}}(q,n) \mid  \prod_{s=1}^m(q-1-i_s)\leq m\}.$$
\end{definition}

\begin{theorem}\label{the3}
Given a positive integer $n$ write $\vec{X}=(X_1, \ldots , X_n)$ and
consider a nonzero polynomial $P(\vec{X}) \in
{\mathbb{F}}_q[\vec{X}]$ with $\deg_{X_i}P<q-1$, $i=1, \ldots ,
n$. Let $m$ be the number of $\vec{x} \in ({\mathbb{F}}_q^\ast)^n$ with
  $P(\vec{x})\neq 0$. Then there does not exist any $(k_1, \ldots
  ,k_n) \in \{0, 1, \ldots, q-2\}^n$ such that 
\begin{eqnarray}
\Supp (X_1^{k_1} \cdots X_n^{k_n} P(\vec{X}) {\mbox{ mod }}
\{X_1^{q-1}-1, \ldots , X_n^{q-1}-1\}) \cap U(q,m,n)=\emptyset. \nonumber
\end{eqnarray}
\end{theorem}

Observe that for $n=1$ we have $U(q,m,n)=\{X^{q-1-m}, \ldots ,
X^{q-1-1}\}$ which is a list of $m$ consecutive monomials. Hence,
Theorem~\ref{the3} is a natural generalization of Theorem~\ref{the1}
to polynomials in more variables.

\noindent {\bf{Proof:}} 
Let $P(\vec{X})$ and $m$ be as in the theorem. Aiming for a
contradiction assume that an $X_1^{k_1} \cdots X_n^{k_n}$ exists such that 
\begin{eqnarray}
\Supp \big( X_1^{k_1} \cdots X_n^{k_n} P(\vec{X}) {\mbox{ mod }}
\{X_1^{q-1}-1, \ldots , X_n^{q-1}-1\} \big) \cap U(q,m,n) =\emptyset . \nonumber
\end{eqnarray}
According to Proposition~\ref{pro2}  
$$X_1^{k_1} \cdots X_n^{k_n} P(\vec{X}) {\mbox{ mod }}
\{X_1^{q-1}-1, \ldots , X_n^{q-1}-1\}$$
has at least $m+1$
non-roots in ${\mathbb{F}}_q^\ast$; and so has $P(\vec{X})$.\\

The generalization of Theorem~\ref{the2} is as follows:

\begin{theorem}\label{the4}
Consider sets $S_i\subseteq {\mathbb{F}}_q^\ast$, $i=1, \ldots ,
n$. Write $s_i=|S_i|$ and assume $0<s_i <q-1$, $i=1, \ldots , n$, not
all $s_i$ being equal to $q-2$. Define
$T_i={\mathbb{F}}_q^\ast\backslash S_i$ and let $t_i=|T_i|$ and
$m=\prod_{i=1}^nt_i$ (by the above assumption on $s_i$ we have $m \geq
2$). Consider 
\begin{equation}
P(\vec{X})=\prod_{i=1}^n\prod_{x \in
  S_i}(X_i-x).\label{eqtrekantm}
\end{equation}
 There exists an $X_1^{k_1}\cdots X_n^{k_n}$ with $0 <
k_1, \ldots , k_n < q-1$ such that
\begin{eqnarray}
\Supp (X_1^{k_1} \cdots
  X_n^{k_n}P(\vec{X}) {\mbox{ mod }} \{X_1^{q-1}-1, \ldots ,
  X_n^{q-1}-1\}) \cap U(q,m-1,n)=\emptyset \label{eqsnabel} 
\end{eqnarray}
if and only if for $i=1, \ldots , n$ it holds that $T_i$ is contained
in $\gamma_iH_i$ for some $\gamma_i\in {\mathbb{F}}_q^\ast$ and for some
proper multiplicative subgroup $H_i$ of ${\mathbb{F}}_q^\ast$.
\end{theorem}

We note that the role of the assumption $m \geq 2$ is to make $U(q,m-1,n)$
non-empty.\\
 
As already observed, for $n=1$ we have $U(q,m-1,n)=\{X_1^{q-m},
\ldots , X_1^{q-2}\}$. Therefore for $n=1$ the assumption~(\ref{eqsnabel}) is equivalent
to saying that 
$${\mathcal{M}}(q,n) \backslash \big( \Supp P \cup
U(q,m-1,n) \big)$$ 
contains a set $D$ such that 
$$U(q,m-1,n) \subseteq
X_1^{k_1}D {\mbox{ mod }} \{X_1^{q-1}-1\}$$ 
(a similar remark does not
hold for $n > 1$.) In other words, for $n=1$, ${\mathcal{M}}(q,n) \backslash
\Supp P$ contains besides $U(q,m-1,n)$ also a translated copy of
$U(q,m-1,n=1)$ which is disjoint from $U(q,m-1,n)$. We have argued that
Theorem~\ref{the4} reduces to Theorem~\ref{the2} in the case that
$n=1$.\\
Turning to the general case of $n \geq 1$ one sees by inspection that $U(q,m-1,n) \subseteq {\mathcal{M}}(q,n)
\backslash \Supp P$. The condition $0<k_i<q-1$, $i=1,\ldots
, n$ means that the sets assumed to exist or proved to exist,
respectively, in Theorem~\ref{the4} are different from $U(q,m-1,n)$
itself; but they may have an overlap with this set.\\

Before giving the proof we illustrate the theorem with an
example.

\begin{example}\label{ex2}
This is a continuation of Example~\ref{ex1} where we considered
polynomials $P(X) \in {\mathbb{F}}_{16}[X]$ of the form~(\ref{eqtrekant}) satisfying the
conditions in Theorem~\ref{the2}. In this example we consider a
polynomial $P(X_1,X_2) \in {\mathbb{F}}_{16}[X_1,X_2]$ of the
form~(\ref{eqtrekantm}) satisfying the condition in
  Theorem~\ref{the4}. Choosing $T_1=\{\alpha^2,\alpha^5,
  \alpha^8,\alpha^{11},\alpha^{14}\}$ and
  $T_2=\{\alpha,\alpha^6,\alpha^{11}\}$ we get that the support of
  $P(\vec{X})$ is 
\begin{eqnarray}
\{1, X_1^5,X_1^{10}, X_2^3, X_1^5X_2^3,X_1^{10}X_2^3,
  X_2^6, X_1^5X_2^6,X_1^{10}X_2^6, \, \, \, \, \, \, \nonumber \\
 X_2^9, X_1^5X_2^9,X_1^{10}X_2^9,
  X_2^{12}, X_1^5X_2^{12},X_1^{10}X_2^{12}\}. \nonumber 
\end{eqnarray}
Clearly, $m=5 \cdot 3=15$. 
In Figure~\ref{figo1} the support is illustrated with diamonds. A
set $D$ is illustrated with filled circles. This set satisfies that $D
\subseteq {\mathcal{M}}(q,n) \backslash \Supp P$ and that  
$$X_1^5X_2^3D {\mbox{ mod }} \{X_1^{q-1}-1, \ldots  X_m^{q-1}-1\}=U(q,m-1,2).$$
Hence, 
$$\Supp (X_1^5X_2^3 P {\mbox{ mod }} \{X_1^{15}-1,X_2^{15}-1\}) \cap
U(q,m-1,2) = \emptyset.$$
\begin{figure}
\begin{center}
\includegraphics[width=50mm]{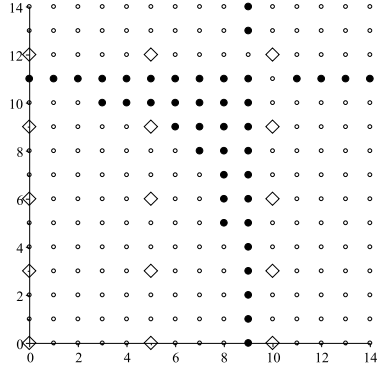}
\end{center}
\caption{The situation in Example~\ref{ex2}.}
\label{figo1}
\end{figure}

\end{example}

It is possible to give a proof of Theorem~\ref{the4} which as a main
tool uses Theorem~\ref{the2} and Proposition~\ref{pro1} in combination with a study of the shape
of $U(q,m-1,n)$. Using this approach the proof of the ``if'' part
becomes straight forward whereas the proof of the ``only if''
part becomes technical and 
requires more care. Instead of stating the technical proof of the
``only if'' part we shall present a self contained
proof of the ``only if'' part based
on the technique from~\cite{kopparty2014roots}. Our proof calls for
the following lemma which has some interest in itself. We state the
lemma in a slightly more general version than shall be needed (we will employ 
the lemma with ${\mathbb{F}}={\mathbb{F}}_q$ and $A_1=\cdots = A_n= {\mathbb{F}}_q^\ast$).

\begin{lemma}\label{lem1}
Given a field ${\mathbb{F}}$ let $A_1, \ldots , A_n \subseteq
{\mathbb{F}}$ be finite sets. Consider proper subsets $B_1 \subsetneq
A_1, \ldots , B_n \subsetneq A_n$ and write
$$P(\vec{X})=\prod_{i=1}^n \prod_{x \in B_i}(X_i-x).$$
Assume that $G(\vec{X}) \in {\mathbb{F}}[\vec{X}]$ is a polynomial
with $\deg_{X_i} G<| A_i|$, $i=1, \ldots , n$ such that
$$\{ \vec{x} \mid \vec{x} \in A_1 \times \cdots \times A_n,
F(\vec{x})=0\} \subseteq \{ \vec{x} \mid \vec{x} \in A_1 \times \cdots \times A_n,
G(\vec{x})=0\}.$$
Then $F(\vec{X})$ divides $G(\vec{X})$.
\end{lemma}

\noindent {\bf{Proof:}} It is enough to prove that $(X_s-x)$ divides $G(\vec{X})$ for
arbitrary $s \in \{1, \ldots , n\}$ and $x \in B_s$. We can write
$G(\vec{X})=Q(\vec{X})(X_s-x)+R(\vec{X})$ where $R(\vec{X})$ is a
polynomial in ${\mathbb{F}}[X_1, \ldots , X_{s-1},X_{s+1}, \ldots ,
  X_n]$ and where $\deg_{X_i} R < |A_i|$ for $i \in \{1, \ldots , s-1,
s+1, \ldots , n\}$. We observe that $(\alpha_1, \ldots ,
\alpha_{s-1},x,\alpha_{s+1}, \ldots , \alpha_n)$ is a root of $P$ and
thereby also of $G$,  for
all $(\alpha_1, \ldots , \alpha_{s-1},\alpha_{s+1}, \ldots ,
\alpha_n)$ in $A_1 \times \cdots \times A_{s-1} \times A_{s+1} \times
\cdots \times A_n$. But then $(\alpha_1, \ldots , \alpha_{s-1},\alpha_{s+1}, \ldots ,
\alpha_n)$ is a root of $R$ and from the Chinese remainder theorem it
follows that $R(\vec{X})=0$.\\

We are now ready to prove Theorem~\ref{the4}.\\

\noindent {\bf{Proof:}} Assume that there
exists an $X_1^{k_1}\cdots X_n^{k_n}$ with $0 <
k_1, \ldots , k_n < q-1$ such that~(\ref{eqsnabel}) holds true.
Write
$$G(\vec{X})=X_1^{k_1} \cdots X_n^{k_n}P(\vec{X}) {\mbox{ mod }}
\{X_1^{q-1}-1, \ldots , X_n^{q-1}-1\}.$$
Clearly the roots of $P(\vec{X})$ in $({\mathbb{F}}_q^\ast)^n$ are
also roots of $G(\vec{X})$. Hence, by Lemma~\ref{lem1}, 
$G(\vec{X})=Q(\vec{X})P(\vec{X})$ for some $Q(\vec{X}) \in
{\mathbb{F}}_q[\vec{X}]$. Recall that $X_1^{s_1}\cdots X_n^{s_n}$ is
the leading monomial of $P(\vec{X})$. If we consider a monomial $N$ such that
$NX_1^{s_1} \cdots X_n^{s_n} \in {\mathcal{M}}(q,n)$ then either $N=1$
or $NX_1^{s_1} \cdots X_n^{s_n} \in U(q,m-1,n)$. From
assumption~(\ref{eqsnabel}) it therefore follows that
$G(\vec{X})=\alpha P(\vec{X})$ for some $\alpha \in
{\mathbb{F}}_q^\ast$. This implies that  
$$P(\vec{X}) (X_1^{k_1} \cdots X_n^{k_n}-\alpha)=0 {\mbox{ mod }}
\{X_1^{q-1}-1, \ldots , X_n^{q-1}-1\}.$$
However, then all non-roots of $P(\vec{X})$ in $({\mathbb{F}}_q^\ast)^n$ --
that is the elements of $T_1 \times \cdots \times T_n$ -- must be
roots of $X_1^{k_1}\cdots X_n^{k_n}-\alpha$. In other words, for
$\vec{x} \in T_1\times \cdots \times T_n$ we have
$x_i^{k_i}=\alpha_i$, $i=1, \ldots , n$ where 
$\prod_{i=1}^n \alpha_i =\alpha$. Consider an $i$ such that $t_i
>1$. Let $y$ and $z$ be two different elements in $T_i$. For fixed
$x_j \in T_j$, $j \in \{1, \ldots , i-1, i+1, \ldots , n\}$ both
$(x_1, \ldots ,x_{i-1},y,x_{i+1}, \ldots , x_n)$ and $(x_1, \ldots
,x_{i-1},z,x_{i+1}, \ldots , x_n)$ satisfy that 
they produce the value $\alpha$
when plugged into
$X_1^{k_1}\cdots X_n^{k_n}$.
Hence, 
$\alpha_1, \ldots , \alpha_n$ are unique. Let $H_i=\{\beta \in
{\mathbb{F}}_q^\ast \mid \beta^{k_i}=1\}$ (which is a proper subgroup
of ${\mathbb{F}}_q^\ast$ as $0 < k_i < q-1$), and $\gamma_i \in
T_i$. Then $T_i \subseteq \gamma_i H_i$.\\
We next prove the ``if'' part of the theorem. Assume that $T_i
\subseteq \gamma_i H_i$, $i=1, \ldots , n$ and write $d_i = |
H_i|$. Define 
$$W_i=\big\{ X_i^v \mid v\in \{jd_i, \ldots , jd_i+(d_i-t_i)\mid j=0,
\ldots , \frac{q-1}{d_i}-1\} \big\}$$
and 
$$W=\{X_1^{v_1} \cdots X_n^{v_n} \mid X_i^{v_i} \in W_i, i=1,
\ldots , n\}.$$
Proposition~\ref{pro1} tells us that $\Supp P \subseteq W$ and by
inspection we find that $U(q-1,m-1,n)$ is contained in
${\mathcal{M}}(q,n)\backslash W$. By symmetry we have
$$X_1^{k_1} \cdots X_n^{k_n} W {\mbox{ mod }} \{X_1^{q-1}-1, \ldots ,
  X_n^{q-1}-1\}=W$$
for all $(k_1, \ldots , k_n)$ where for $i=1, \ldots , n$,
$k_i=\ell_id_i$ for some $\ell_i$. The theorem follows.

\section*{Acknowledgements}
This work was supported by the Danish Council for Independent
Research, grant no.\ DFF-4002-00367.

\bibliographystyle{plain}

\end{document}